\documentclass[runningheads]{llncs}
\usepackage{graphicx}
\usepackage{amsmath}
\usepackage{amssymb}
\usepackage[inline]{enumitem}
\usepackage[mode=buildnew]{standalone}%
\usepackage{url}
\newcommand{\mZ}{\mathcal{Z}}\newcommand{\mX}{\mathcal{X}}
\newcommand{\mE}{\mathcal{E}}
\newcommand{\mR}{\mathcal{R}}

\newcommand{\mO}{\mathcal{O}}
\newcommand{\mS}{\mathcal{S}}
\newcommand{\mY}{\mathcal{Y}}

\newcommand{\mb}{\mathbf}

\usepackage{todonotes}

\usepackage[T1]{fontenc}
 \usepackage{tikz}
\usetikzlibrary{arrows.meta}
\usepackage[scaled=0.90]{FiraMono}

\PassOptionsToPackage{hyphens}{url}\usepackage{hyperref}

\usepackage{xcolor}
\hypersetup{
    colorlinks,
    linkcolor={red!50!black},
    citecolor={blue!50!black},
    urlcolor={blue!80!black}
}
\usepackage[noabbrev,capitalize,nameinlink]{cleveref}[0.21]

\usepackage{subcaption}

\begin{document}
\title{Network Design for the Traffic Assignment Problem with Mixed-Integer Frank-Wolfe}
\titlerunning{Network Design for Traffic Assignment with Mixed-Integer Frank-Wolfe}
%
\author{Kartikey Sharma\inst{1}\orcidID{0000-0001-6736-4827} \and
Deborah Hendrych\inst{1,2}
\orcidID{0000-0003-0705-1356} \and
Mathieu Besançon\inst{1,3}
\orcidID{0000-0002-6284-3033} \and
Sebastian Pokutta\inst{1,2}
\orcidID{0000-0001-7365-3000}}
\authorrunning{K.~Sharma et al.}
%
\institute{Zuse Institute Berlin, Takustr. 7, 14195 Berlin, Germany \\
\email{\{kartikey.sharma,hendrych,besancon,pokutta\}@zib.de} \and
Technische Universit\"at Berlin, Straße des 17. Juni 135, 10623 Berlin, Germany \and
Université Grenoble Alpes, Inria, LIG, Grenoble, France
}

\maketitle              
\begin{abstract}
We tackle the network design problem for centralized traffic assignment, which can be cast as a mixed-integer convex optimization (MICO) problem.
For this task, we propose different formulations and solution methods
in both a deterministic and a stochastic setting in which the demand is unknown in the design phase.
We leverage the recently proposed Boscia framework, which can solve MICO problems when the main nonlinearity stems from a differentiable objective function.
Boscia tackles these problems by branch-and-bound with continuous relaxations solved approximately with Frank-Wolfe algorithms.

We compare different linear relaxations and the corresponding subproblems solved by Frank-Wolfe, and alternative problem formulations to identify the situations in which each performs best.
Our experiments evaluate the different approaches on instances from the Transportation Networks library and highlight the suitability of the mixed-integer Frank-Wolfe algorithm for this problem. In particular, we find that the Boscia framework is particularly applicable to this problem and that a mixed-integer linear Frank-Wolfe subproblem performs well for the deterministic case, while a penalty-based approach, with decoupled feasible regions for the design and flow variables, dominates other approaches for stochastic instances with many scenarios.


\end{abstract}
\section{Introduction}

Given a network $G = (V,\mE)$ and a collection of flows between source $\mO \subseteq V$ and destination $\mZ \subseteq V$ vertices  in the network, the traffic assignment problem (TA) finds the equilibrium traffic flow which minimizes the average traffic flow cost and can be expressed as:
\begin{align}
    \min_{\mathbf{x}} &\; c(\mb{x}) := \sum_{e \in \mE} c_e(x_e) \tag{TA}\label{eq:ta}\\
    \text{s.t.} &\; x_e = \sum_{z \in \mZ}x_e^z \;&&\forall e \in \mE\nonumber \\
    &\; \mathbf{x}^z \in \mX^z = \begin{cases}
        \sum_{e \in \delta^{+}(i)}x_e^z - \sum_{e \in \delta^{-}(i)}x_e^z = 0, \;\forall i \in V \ \backslash \ (\mO \cup \mZ)\\
        \sum_{e \in \delta^{+}(i)}x_e^z = d_i^z \;\forall i \in \mO\\
        \sum_{e \in \delta^{-}(z)}x_e^z = \sum_{i \in \mO} d_i^z\\
        \end{cases}\;&&\hspace{-1mm}\forall z \in \mZ.\nonumber
\end{align}

Here, the variable $\mathbf{x}$ captures the flow on the edges of the network with $x_e^z$ being the flow on edge $e$ towards destination $z$ and $d_i^z$ being the expected traffic or traffic demand from source $i \in \mO$ to destination $z \in \mZ$.
The set $\mX^z$ consists of the network flow constraints for destination $z$.
The objective term $c_e(x_e)$ is the traffic cost on arc $e$ which estimates the travel time and is of the form:
\begin{align*}
c_e(x_e) = \alpha_e + \beta_e x_e + \gamma_e x_e^{\rho_e},
\end{align*}
where $\alpha_e, \beta_e \text{ and } \gamma_e$ are constants. The cost increases superlinearly as a function of the total flow $x_e$ on the arc $e$ with power $\rho_e > 1$ thus capturing the congestion effect.
The first constraint aggregates the total flow on the arc to all destinations. 
The set $\mX^z$ captures flow conservation at the origin and intermediate nodes and at destination $z$. 
Here, $\delta^{+}(i)$ and $\delta^{-}(i)$ are the set of outgoing and incoming edges at node $i$ respectively. 
We use $F$ to denote the feasible region of Problem~\eqref{eq:ta}.


In the presence of congestion, the cost of the flow on an arc increases superlinearly as a function of the cumulative flow on the arc. 
At equilibrium, every path between a source-destination pair through which a flow passes must have the same cost, and no other path must be better~\cite{leblanc1975algorithm}.
As such, the traffic flow problem is commonly used as a subproblem to model transportation and communication networks to simulate complex dynamics~\cite{hackl2019estimation}.

One well-established way of tackling Problem~\eqref{eq:ta} is through Frank-Wolfe (FW) or Conditional Gradient algorithms \cite{frank1956algorithm,jaggi2013revisiting,levitin1966constrained},
see for instance \cite{chen2002faster,fukushima1984modified,lee2001accelerating,mitradjieva2013stiff}.
Frank-Wolfe is a first-order method for constrained optimization which iteratively computes a gradient of the nonlinear objective function and solves a linear problem over the constraint set, using the gradient at the current iterate as a linear objective.
For the TA application, one core interest is that the linear subproblem can be solved efficiently by calls to a specialized shortest-path algorithm instead of a generic linear programming solver.
We will exploit this property in particular for some of our solution methods.
We refer the reader to \cite{bomze2021frank,braun2022conditional} for recent overviews of FW methods.

\subsubsection{Network Design for Traffic Assignment.}

\begin{figure}[ht]
    \centering
\begin{subfigure}{.4\textwidth}
  \centering
      \begin{tikzpicture}[>=Stealth, node distance=2cm, font=\Large]

        \node[rectangle, draw, fill=white, text=red] (s1) at (0,2) {$S_1$};
        \node[circle, draw, fill=white, right=of s1] (n1) {1};
        \node[circle, draw, fill=white, below=of n1] (n2) {2};
        \node[circle, draw, fill=white, right=of n1] (n3) {3};
        \node[circle, draw, fill=white, right=of n3] (n4) {4};
        \node[circle, draw, fill=white, below=of n4] (n5) {5};
        \node[circle, draw, fill=white, below=of n2, text=blue] (d) {$D$};
        \node[rectangle, draw, fill=white, below=of n3, text=red] (s2) {$S_2$};

        \draw[->, line width=1pt] (s1) -- (n1) node[midway, above] {$1$};
        \draw[->, line width=0.5pt, dashed] (n1) -- (n2) node[midway, left] {$0$};
        \draw[->, line width=1pt] (n1) -- (n3) node[midway, above] {$1$};
        \draw[->, line width=2.5pt] (n3) -- (n4) node[midway, above] {$2$};
        \draw[->, line width=2.5pt] (n4) -- (n5) node[midway, left] {$2$};
        \draw[->, line width=2.5pt] (n5) -- (d) node[pos=0.4, below=0.15cm] {$2$};
        \draw[->, line width=0.5pt] (n2) -- (d) node[midway, left] {$0$};
        \draw[->, line width=1pt] (s2) -- (n3) node[midway, left] {$1$};

    \end{tikzpicture}
  \caption{Initial network}
  \label{fig:AOptimalTermination}
\end{subfigure}%
\begin{subfigure}{.4\textwidth}
  \centering
      \begin{tikzpicture}[>=Stealth, node distance=2cm, font=\Large]

        \node[rectangle, draw, fill=white, text=red] (s1) at (0,2) {$S_1$};
        \node[circle, draw, fill=white, right=of s1] (n1) {1};
        \node[circle, draw, fill=white, below=of n1] (n2) {2};
        \node[circle, draw, fill=white, right=of n1] (n3) {3};
        \node[circle, draw, fill=white, right=of n3] (n4) {4};
        \node[circle, draw, fill=white, below=of n4] (n5) {5};
        \node[circle, draw, fill=white, below=of n2, text=blue] (d) {$D$};
        \node[rectangle, draw, fill=white, below=of n3, text=red] (s2) {$S_2$};

        \draw[->, line width=1pt] (s1) -- (n1) node[midway, above] {$1$};
        \draw[->, line width=1pt] (n1) -- (n2) node[midway, left] {$1$};
        \draw[->, line width=0.5pt] (n1) -- (n3) node[midway, above] {$0$};
        \draw[->, line width=1pt] (n3) -- (n4) node[midway, above] {$1$};
        \draw[->, line width=1pt] (n4) -- (n5) node[midway, left] {$1$};
        \draw[->, line width=1pt] (n5) -- (d) node[pos=0.4, below=0.15cm] {$1$};
        \draw[->, line width=1pt] (n2) -- (d) node[midway, left] {$1$};
        \draw[->, line width=1pt] (s2) -- (n3) node[midway, left] {$1$};

    \end{tikzpicture}
  \caption{Network with the new edge}
  \label{fig:DOptimalTermination}
\end{subfigure}
\caption{Example of a network with an optional purchasable (dashed) edge. Each edge is labeled with the corresponding flow, there is a unit flow from each source $S_{1,2}$ to the destination $D$.
Because of the superlinear congestion costs, edges with two units of flow incur a cost greater than twice that of a single unit.}
\label{fig:TermiantionOverTimeOptimal}
\end{figure}
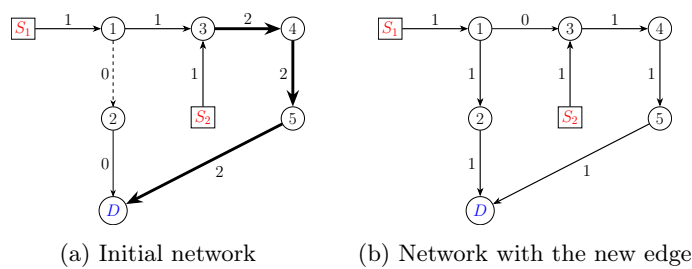

Network design is a generic combinatorial problem that consists of 
choosing a set of edges to add to a network in order to minimize the weighted sum of
a design cost from the additional edges~\cite{johnson1978complexity}
and an operational cost on the resulting network.
We consider in particular the network design problem where the operational cost corresponds to the traffic assignment problem. 
The goal is to determine a set of new edges to add to a network and to assign a set of flows on the resulting graph such that the sum of the traffic flow cost $c(\mathbf{x})$ and the network design cost $\mb{r}^{\top}\mb{y}$ is minimized:
\begin{align}
    \min_{\mathbf{y}, \mathbf{x}} \;\; & \mathbf{r}^\intercal \mathbf{y} + c(\mathbf{x}) \tag{ND}\label{eq:ndta}\\
    \text{s.t.} \;\; & y_e = 0 \Rightarrow x_e \leq 0 && \forall e \in \mR\ \nonumber\\
    & \mb{x} \in F\nonumber \\
    & \mathbf{y} \in \mathcal{Y} \subseteq \{0,1\}^{|\mathcal{R}|}.\nonumber
\end{align}
Here, the subset $\mR \subseteq \mE$ is the collection of edges that can be added to the network and operated on.
An illustration of this problem is given in \cref{fig:TermiantionOverTimeOptimal}.

\noindent
Note that our model is an approximation of the traditional \emph{user} network design problem which separates the $\mb{x}$ and $\mb{y}$ decisions into distinct stages and models~\eqref{eq:ndta} as a bilevel optimization problem. 
The first stage decision chooses the optimal network design $\mb{y}$ to minimize the traffic congestion. 
However, the traffic assignment decision $\mb{x}$ instead of being centrally assigned is calculated by solving the Problem~\eqref{eq:ta} in a second stage given the choice of $\mb{y}$. 
The key distinction in our model is that we simultaneously choose both $\mb{y}$ and $\mb{x}$ whereas the \emph{user}-based approach only consists of choosing $\mb{y}$.
In general, our approach can be used to develop solutions that are an approximation or an initial point for true user equilibrium solutions~\cite{wang2013global} or to lower bound the user network design problem in a Branch \& Bound approach~\cite{leblanc1975algorithm,farvaresh2013branch}.
Our problem would also  be applicable in settings where the flow decisions are centralized such as freight or rail~\cite{feng2024integrating,uddin2015freight}.


\subsubsection{Stochastic Network Design.}
In addition to deterministic network design, in which the demands are known beforehand, we consider a stochastic variant with uncertain demands.
In this stochastic network design problem, the decision-maker designs the network, after which the demand is revealed and the optimal flow can be computed.
We can write this problem as:
\begin{align}
    \min_{\mathbf{y}, \mathbf{x}} \;\;& \mathbf{r}^\intercal \mathbf{y} + \sum_{s \in \mS} p_s c(\mathbf{x}_s) \tag{SND}\label{eq:sndta}\\
    \text{s.t.} \; & y_e = 0 \Rightarrow x_{e,s} \leq 0 && \forall e \in \mR,\; \forall s \in \mS\nonumber\\
    & \mb{x}_s \in F_s\; &&\forall s \in \mS\nonumber\\
    & \mathbf{y} \in \mathcal{Y} \subseteq \{0,1\}^{|\mathcal{R}|}.\nonumber
\end{align}
In the above formulation, scenario $s \in \mS$ occurs with probability $p_s$. 
There are separate flow variables $\mb{x}_s$ for each scenario. 
The network flow feasible region $F_s$ is scenario-dependent due to the presence of uncertain demand.
The network design decisions and constraints are the only elements linking the problem across the scenarios, a property we exploit in solution methods.

\subsubsection{Previous Work.}
Traditional network design problems have primarily focused on user equilibrium flows in response to network design decisions. 
The resulting problem is modeled as an integer bilevel optimization problem and solved using approaches such as Branch-and-Bound~\cite{farvaresh2013branch,leblanc1975algorithm}, Benders' Decomposition~\cite{fontaine2014benders,luathep2011global} and Lagrangian Relaxation-based methods~\cite{gendron2019revisiting}. 
Other approaches have focused on reformulating the bilevel problem as a single-stage Problem~\cite{farvaresh2011single} or leveraging heuristics such as Simulated Annealing and Swarm Optimization~\cite{poorzahedy2007hybrid}.
Recent work has also focused on leveraging global optimization algorithms to solve network design problems~\cite{wang2015novel,wang2013global}.
The centralized decision-making process in our problem is more similar to a fixed-charge network flow (FCNF) problem~\cite{bruynooghe1972optimal,hirsch1968fixed} but one where the fixed costs only apply to a subset of the arcs. 
Similar to the traditional network design, FCNF problems have been solved using iterative algorithms such as Branch-and-Bound~\cite{ortega2003branch} and Benders decomposition~\cite{fragkos2021decomposition}.  
Due to the large size of the full-scale FCNF problem, a variety of heuristics have also been developed to obtain solutions~\cite{hewitt2010combining,xie2012nonlinear}.

The problem we consider straddles the line between these two settings as we combine the nonlinear objective of traditional network design with the fixed costs of the FCNF problem. 
Existing research on this topic has focused on either quadratic~\cite{xie2012nonlinear} or piecewise linear costs~\cite{gendron2017reformulations}.  




\subsubsection{Contributions.}
We construct several formulations of the centralized network design problem for traffic assignment \eqref{eq:ndta} and its stochastic variant that selects the best additional edges while minimizing the expected operational cost over a finite set of scenarios~\eqref{eq:sndta}. 
Specifically, we leverage 
\begin{enumerate*}[label=(\roman*)]
\item a mixed-integer conic formulation solved with an off-the-shelf solver, 
\item a mixed-integer formulation that maintains the original problem structure over a mixed-integer feasible set, 
\item and a penalty-based formulation that decomposes the mixed-integer set, relaxing the constraints and adding a nonlinear penalty in the objective.
\end{enumerate*}

We then design solution methods tailored to the different formulations.
In particular, we base multiple methods on the mixed-integer Frank-Wolfe algorithm based on the recently proposed \texttt{Boscia.jl} \cite{hendrych2023convex} framework.
This allows us to exploit the combinatorial structure of the linearized subproblems in the mixed-integer nonlinear formulations.
Two of these approaches solve the Mixed-Integer Programming (MIP) subproblem using direct and cut generation methods. 
We also propose an alternative approach that replaces the design constraint with a nonlinear penalty on its violation.
This splits the feasible set into the design variable set $\mathcal{Y}$ and the flow polytope, both of which can benefit from tailored algorithms to solve the corresponding linear subproblems.
Along with the mixed-integer Frank-Wolfe approach, we also attempt to solve the MICO problem directly using conic and nonlinear optimization approaches.
Note that for all of the Frank-Wolfe-based approaches we use a variant called the \emph{Blended Pairwise Conditional Gradient} (BPCG) which performs better.

Finally, we assess the different solution methods on deterministic and stochastic traffic assignment instances.
The results show that the MIP-based methods outperform other methods for deterministic instances and small numbers of scenarios while the penalty-based method is better for more scenarios.

\section{Solution Methods}
This section discusses the different formulations and corresponding solution methods that we use to solve the network design problem. 
Three of these methods leverage the Boscia framework and consist of different linear minimization oracles (LMOs) whereas the last two directly solve the optimization problem. 

\subsection{Integer Frank-Wolfe}
The Integer Frank-Wolfe (IFW) approach uses off-the-shelf Mixed-Integer Linear Programming (MILP) solvers to solve a linearized version of Problem~\eqref{eq:ndta}.
This linear problem computes the descent direction $\mb{v}_y, \mb{v}_x$ at any given iteration $t$, and in it, the original nonlinear objective has been replaced by a linear objective using the gradient at the current iterate $\mb{x}_t$. 
The optimal solution to this problem is the best feasible descent direction. 
We then take a step using this direction.
\begin{align}
    \min_{\mathbf{v}_y, \mathbf{v}_x} &\; \mathbf{r}^\intercal \mathbf{v}_y + \nabla c(\mathbf{x}_t)^{\top} \mb{v}_x \tag{IFW}\label{eq:lndta}\\
    \text{s.t.} &\; v_{e,x}^z \leq M^z v_{e,y}&& \forall e \in \mR,\; \forall z \in \mZ \nonumber\\
    &\; \mb{v}_x \in F \nonumber\\
    &\; \mathbf{v}_y \in \mathcal{Y} \subseteq \{0,1\}^{|\mathcal{R}|}.\nonumber
\end{align}

The standard FW method is known to suffer slow convergence rates over polytopes, even when the objective function is strongly convex. 
Furthermore, it requires one LMO call per iteration, which is prohibitively expensive in many settings, including ours, where the LMO requires solving a combinatorial problem (min-cost flow for the penalty approach or worse, linear network design for IFW).
Therefore, we use the BPCG variant \cite{tsuji2022pairwise}, with two advantageous effects.
The first is improved convergence by choosing the best direction between a FW step and a pairwise step which works directly on the weights of the convex combination forming the current iterate.
The second being~\emph{lazification}, i.e.,~the ability of the algorithm to select vertices that are not solutions to the LMO at the current iterate as progress directions, as long as they provide sufficient progress.

\subsection{Penalty-based Integer Frank-Wolfe}
A key difficulty in solving the Integer LMO and hence Problem~\eqref{eq:ndta} in an efficient fashion is the presence of the network design constraint.
Without this constraint, the problem is separable for each destination $z \in \mZ$ (and each scenario $s \in \mS$).
As such, we relax this constraint and instead penalize its violation.
Furthermore, since we are using an inherently nonlinear approach, the penalty term can also be nonlinear.
The resulting problem can be expressed as:
\begin{align}
    \min_{\mathbf{y}, \mathbf{x}} &\; \mathbf{r}^\intercal \mathbf{y} + c(\mathbf{x}) + \mu \sum_{z \in \mZ} \sum_{e \in \mR} \max(x_e^z - M^z y_e, 0)^p  \tag{PB-ND}\label{eq:rndta}\\
    \text{s.t.} 
    &\; \mb{x} \in F\nonumber\\
    &\; \mathbf{y} \in \mathcal{Y} \subseteq \{0,1\}^{|\mathcal{R}|}.\nonumber
\end{align}

In Problem~\eqref{eq:rndta}, we have moved the constraint which toggles the arcs into the objective and have penalized its violation the extent of which is parameterized by $\mu$ and $p$.
The key benefit of this relaxation can be seen when we write down the linear minimization oracle for this problem:
\begin{align}
    \min_{\mathbf{v}_y, \mathbf{v}_x} &\; (\mathbf{r} + \mb{g})^\intercal \mathbf{v}_y + (\nabla c(\mathbf{x}_t) + \mb{h})^\intercal \mb{v}_x   \tag{PB}\label{eq:lrndta}\\
    \text{s.t.} 
    &\; v_{e,x} = \sum_{z \in \mZ}v_e^z && \forall e \in \mE\nonumber \\
    &\; \mathbf{v}^z \in \mX^z \;&&\forall z \in \mZ\nonumber\\
    &\; \mathbf{v}_y \in \mathcal{Y} \subseteq \{0,1\}^{|\mathcal{R}|}.\nonumber
\end{align}
where 
\begin{align*}
\mb{g} &= \{g_e\}_{e \in \mR} := - 2 \mu p \sum_{z \in \mZ} M^z \max(x_{e,t}^z - M^z y_{e,t},0)^{p-1} \\
\mb{h} &= \{h_e^z\}_{e \in \mE, z \in \mZ} := 
\begin{cases}
2 \mu p \max(x_{e, t}^z - M^z y_{e, t},0)^{p-1} & \text{ if } e \in \mR\\
0 & \text{ if } e \notin \mR.
\end{cases}
\end{align*}

With this setup, we can separate the Problem~\eqref{eq:lrndta} into two separate subproblems over $\mb{v}_x$ and $\mb{v}_y$ which can be solved much more efficiently.
For the problem over $\mb{v}_x$, we can utilize algorithms for the \emph{Shortest Path Problem}. 
The problem over $\mb{v}_y$ can be solved with an appropriate method depending on $\mY$.
This allows for very fast computation. 
However, this approach has two key limitations. First, the parameter $\mu$ needs to be estimated to solve the problem correctly. 
Second, the relaxation of the network design constraints means that we lose structural information as the new descent direction for $\mb{x}$ i.e., $\mb{v}_x$ only depends upon the previous value of $\mb{y}_t$ and not the current value $\mb{v}_y$.
As such, solving the overall problem with this LMO requires more iterations. 

\subsection{Benders decomposition of the linear network design}

An essential part of the computations for IFW consists of solving the linearized subproblems whose structure is akin to linear network design.
An established approach for these problems is Benders decomposition, removing the continuous flow variables and solving for the design variables $\mathbf{v}_y$ with an exponential number of inequalities separated on the fly, we refer the reader to \cite{rahmaniani2017benders} for a review of Benders decomposition techniques and applications.
In particular, this decomposition allows us to leverage iterative calls to network flow subproblems to generate valid inequalities.
For this approach, we iterate between the master problem containing the design and auxiliary epigraph variables, and subproblems.
These can be abstractly represented as:

\vspace{-5mm}
\begin{minipage}[t]{0.4\textwidth}
\begin{align*}
    \max_{\mathbf{r}, \mathbf{t}, \mb{s}, \mb{p}} &\; (\mb{r} - \mb{t})^{\intercal} \mb{d} - \mb{s}^{\intercal} \mb{v}_y \nonumber\\
    \text{s.t.} &\; (\mathbf{r}, \mathbf{t}, \mb{s}, \mb{p}) \in \mathcal{D}, \nonumber
\end{align*}
\end{minipage}
\begin{minipage}[t]{0.4\textwidth}
\begin{align*}
    \min_{\mathbf{v}_y, \eta} &\; \mathbf{r}^\intercal \mathbf{v}_{y} + \eta \qquad\tag{BDM} \\
    \text{s.t.} &\; 
    \eta + \mb{s}_k^{\intercal} \mb{v}_y \geq (\mb{r}_k - \mb{t}_k)^{\intercal} \mb{d}\; &&\forall k \in O  \nonumber \\
    &\; \mb{s}_k^{\intercal} \mb{v}_y \geq (\mb{r}_k - \mb{t}_k)^{\intercal} \mb{d} \;&&\forall k \in I \nonumber \\
    &\; \mathbf{v}_y \in \mathcal{Y} \subseteq \{0,1\}^{|\mathcal{R}|}.\nonumber
\end{align*}
\end{minipage}
\\

\noindent
Here, the sets $O$ and $I$ are a collection of optimality and feasibility cuts and $\mathcal{D}$ is the dual feasible region for the subproblem.
The above formulation offers the simplest setup of Benders' decomposition and can be improved by leveraging the various techniques that have been developed for such problems over the years (e.g., stabilization). We leave these aspects to future work.

\subsection{Direct Solution Methods}
Along with the previously discussed FW-based approaches, we also attempt to directly solve problem~\eqref{eq:ndta} using a mixed-integer non-linear solver. 
A line of work that has proven effective for some mixed-integer nonlinear problems is to reformulate the nonlinearities using conic constraints \cite{aps2020mosek,coey2022solving,coey2020outer}.
We therefore build mixed-integer conic formulations using power cone constraints to reformulate the nonlinear objective of~\eqref{eq:ndta}.
Let $\mathcal{K}_{\omega}$ denote the power cone of exponent $\omega$:
\begin{align*}
    \mathcal{K}_{\omega} = \{(a,b,c) \in \mathbb{R}^3 \,|\, a^{\omega} b^{1-\omega} \geq |c| \}.
\end{align*}
The problem can then be represented as:
\begin{align*}
    \min_{\mathbf{y}, \mathbf{w}, \mathbf{x}} &\; \mathbf{r}^\intercal \mathbf{y} + \sum_{e\in\mE } \alpha_e + \beta_e x_e + \gamma_e w_e \\
    \text{s.t.} &\; y_e = 0 \Rightarrow x_e \leq 0 \;&& \forall e \in \mR\\
    &\; \mb{x} \in F\\
    &\; (w_e, 1, x_e) \in \mathcal{K}_{1/\rho_e} && \forall e \in \mE \\
    &\; \mathbf{y} \in \mathcal{Y} \subseteq \{0,1\}^{|\mathcal{R}|}.
\end{align*}
The conic constraints are derived from an epigraph of the nonlinear objective:
\begin{align*}
& w_{e} \geq x_{e}^{\rho_e}
\Leftrightarrow \;\; w_{e}^{1/\rho_e} \geq x_{e}
\Leftrightarrow \;\; (w_e, 1, x_{e}) \in \mathcal{K}_{1/\rho_e} \quad \forall e \in \mE.
\end{align*}
The reformulated problem can be tackled by mixed-integer conic solvers such as \texttt{Pajarito.jl}~\cite{coey2020outer}.
Furthermore, the formulation can be improved with a perspective reformulation representing the constraint $y_{e} = 0 \Rightarrow x_{e} = 0$ jointly with the nonlinear function corresponding to the contribution of $x_e$ to the objective if $y_e = 1$, see \cite{frangioni2006perspective,lee2022gaining} for the fundamentals of perspective functions and their power cone reformulations.
This results in the following conic formulation:
\begin{align*}
\min_{\mathbf{y},\mathbf{w}, \mathbf{x}}\;\; & \mathbf{r}^\intercal \mathbf{y} + \sum_{e\in\mathcal{E}} \alpha_e + \beta_e x_e + \gamma_e w_e \\
\text{ s.t. } & (w_{e}, y_{e}, x_{e}) \in \mathcal{K}_{1/\rho_e} && \forall e \in \mR \\
& (w_{e}, 1, x_{e}) \in \mathcal{K}_{1/\rho_e} && \forall e \in \mathcal{E} \ \backslash \ \mR \\
& \mathbf{x} \in F, \mathbf{y} \in \mathcal{Y}.
\end{align*}

Besides conic optimization, we also attempted to solve problem~\eqref{eq:ndta} directly using the SCIP solver which can solve mixed-integer nonlinear problems.


\section{Computational Experiments}
In this section, we compare the performances of Boscia for the different linear minimization oracles, the conic formulation, and the mixed-integer nonlinear formulation. 
To create the necessary instances we used datasets from the the Transportation Networks library~\cite{transportnetlibrary}.
We generated the network design problem by selecting a fraction of the arcs and removing them. 
The price to add them back was set to the equilibrium traffic flow cost divided by the number of edges in the network. 
We created scenarios for the stochastic model by multiplying the nominal (deterministic) demand between any source and destination pair by a random number in the interval $[1.0, 1.1]$. 
We make the code for these experiments available on GitHub\footnote{\url{https://www.github.com/kartikeyrinwa/Network\_Design\_with\_Integer\_Frank\_Wolfe}}.

\subsubsection{Instances.}
We consider 5 datasets from the Transportation Networks library: \emph{Berlin-Friedrichshain}, \emph{Berlin-Tiergarten}, \emph{Berlin-Prenzlauerberg-Center}, \emph{Berlin-Mitte-Center} and \emph{Anaheim}. 
These networks all have fewer than 1000 arcs and are solvable for stochastic settings. 
The network design instances are created by removing $1-5\%$ of the arcs, leading to $25$ instances. 
The stochastic instances are constructed with $2$ to $50$ scenarios from the deterministic instance. 
We consider an instance as solved if the duality gap is less than $5\%$ and the constraint violation (for the penalty method) is less than $0.01$.
The constraint violation is calculated as the maximum flow over all removed arcs and across all scenarios. 

\subsubsection{Stochastic Network Design.}
For the IFW approach, the stochastic problem is a linearized version of Problem~\eqref{eq:sndta}.
For the PB and BDM approaches the linearized problem is solved separately for each scenario given the design decisions.
We use a penalty \(\mu = 1000.0\) and power \(p = 1.5\) for this approach.
The Conic and MINLP approaches solve Problem~\eqref{eq:sndta} directly.

\subsubsection{Results.}
We compare the Integer Frank-Wolfe (IFW) approach, using both SCIP~\cite{scip_solver} and HiGHS~\cite{huangfu2018parallelizing} to solve the MIP subproblems, the penalty-based method (PB), and the Benders Decomposition LMO (BDM). 

\begin{figure}
    \centering\includegraphics[width=4.5cm]{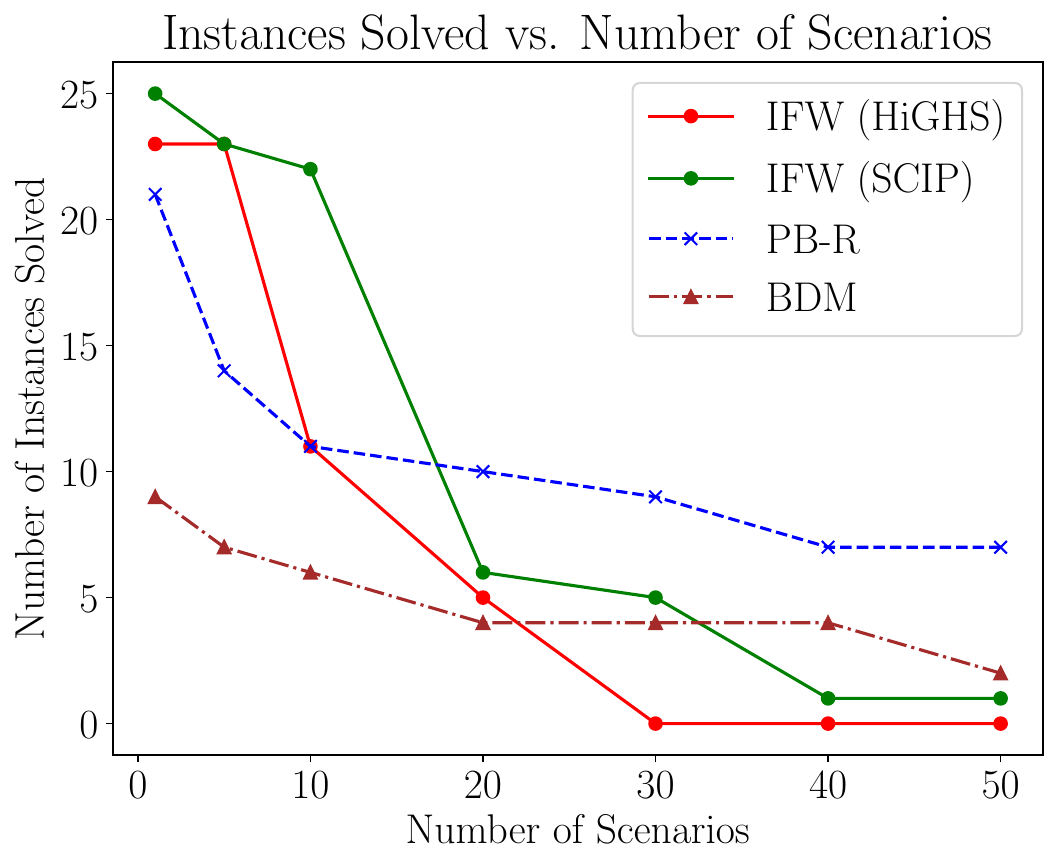}
    \includegraphics[width=4.5cm]{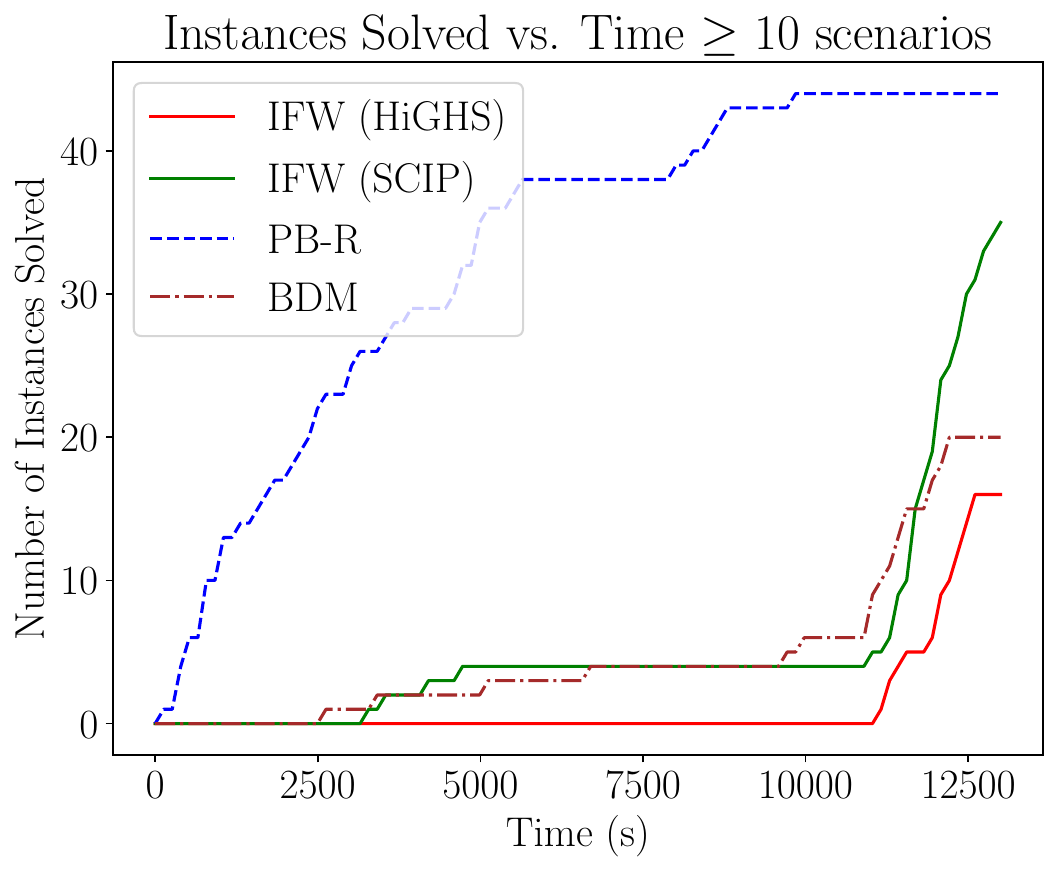}
    \caption{Comparison of the Integer LMO, the Shortest path LMO, and the Benders Decomposition based LMO with solutions returned within 13000s.}
    \label{fig:bndlmo_late}
\end{figure}

\cref{fig:bndlmo_late} (left) shows that for the deterministic problem and instances with a few scenarios, IFW outperforms all of the other methods in terms of instances solved. The performance of the method degrades as the number of scenarios increases; this holds for both MIP solvers.
The penalty-based method outperforms all others for instances with a large number of scenarios, corresponding to a large number of flow variables. Compared to IFW, the penalty-based method results in inexpensive LMO calls, but looser continuous relaxations, yielding large branching trees whereas the IFW approach has much smaller trees due to its tighter formulation.
The Benders' decomposition performs worse than the other methods on instances with few scenarios but improves over the IFW method for larger scenario counts.
\cref{fig:bndlmo_late} (right) depicts the number of instances solved as a function of time all instances with more than 10 scenarios are aggregated.
The plot shows that for instances with more scenarios, the penalty-based method significantly outperforms the other methods, solving many more instances.
We also observe that BDM performs better than the direct MIP formulation for higher scenario counts, but is still outperformed by the penalty-based method in this regime, suggesting that the strength gained in the relaxation is not sufficient to balance out its computational cost.

We also compared the different ways of modeling the network design constraint in Problem~\eqref{eq:ndta}. 
This was a performance comparison between the basic network design constraints, indicator constraints, and both the network design and indicator constraints.
All three methods performed similarly to each other. 

\paragraph{Mixed-Integer Conic Formulation.}
We also tackled network design with the standard and perspective conic formulations.  
For this, we used \texttt{Pajarito.jl} which implements a polyhedral outer approximation with cuts generated from the dual conic constraints using~\texttt{Hypatia.jl}~\cite{coey2022solving}. 
The conic methods were unable to solve even the smallest instances used in the experiments so far.
As such, we created a smaller example from the \emph{Berlin-Friedrischain} dataset by reducing the number of source-destination flows to 6 out of the original 506 and removing 40\% of the arcs. 
We then considered $1\%$ removable arcs with $1,2,$ and $5$ scenarios.

\begin{table}
\centering
\subfloat[Standard]{
\begin{tabular}{|l|l|l|l|}
\hline
Scenarios         & 1    & 2    & 5    \\ \hline
Time (s) & 90   & 305  & 3078 \\ \hline
Cuts     & 1956 & 3910 & 9796 \\ \hline
\end{tabular}}
\quad
\subfloat[Perspective]{
\begin{tabular}{|l|l|l|l|}
\hline
Scenarios         & 1    & 2    & 5    \\ \hline
Time (s) & 92   & 305  & 3688 \\ \hline
Cuts     & 1960 & 3914 & 9787 \\ \hline
\end{tabular}}
\caption{Time taken and cuts generated by the conic formulations.}
\end{table}
Table~1 shows the time taken and the cuts required to solve the standard and perspective conic formulations on the toy example. 
These approaches, even in a simplified example, are much slower than the other methods.
This shows the challenges for the outer approximation method with conic constraints.

\paragraph{Mixed-Integer Nonlinear Formulation.}
Finally, we attempted to solve the network design problems directly with SCIP used as an MINLP solver (which requires the symbolic expression of the nonlinear objective).
The resulting setup successfully solved toy instances (such as those solved by the conic formulation) in a few seconds. 
However, the solver was much slower than our approaches for practical instances and also suffered significant numerical issues in the LP solves.

\section{Conclusion}

In this paper, we developed new solution methods to solve network design problems for centralized traffic assignment based on modern algorithmic frameworks, namely Branch-and-Bound with Frank-Wolfe and mixed-integer conic optimization.
For the former, we developed three problem formulations.
The first two use an MILP to capture the network design and flow decisions and leverage either an MIP solver or cut generation to solve the problem.
The last expresses the network design constraints through a nonlinear penalty function, making the linear subproblems fully separable and exploiting specialized combinatorial algorithms.

Our computational experiments show that the MIP-based Frank-Wolfe formulation outperforms other approaches for deterministic instances and stochastic ones with fewer scenarios, benefitting from a tight convex relaxation and resulting in small branching trees.
On stochastic instances with more scenarios, the penalty-based formulation showcases better performance, solving or starting to close the gap for the hardest instances on which even a single LMO call of IFW is challenging.
The LMO based on the Benders decomposition approach, while not as good as the previous two, also performed well and has significant potential for improvement using better cuts. 

\subsubsection{\ackname} Research reported in this paper was partially supported through the Research Campus Modal funded by the German Federal Ministry of Education and Research (fund numbers 05M14ZAM,05M20ZBM) and the Deutsche Forschungsgemeinschaft (DFG) through the DFG Cluster of Excellence MATH+.

%
%
%
\bibliographystyle{splncs04}
\bibliography{references}
\end{document}